\documentclass[12pt]{amsart}
\advance\hoffset by -0.5 truecm
\usepackage{amssymb}
\newtheorem{Theorem}{Theorem}[section]
\newtheorem{Lemma}[Theorem]{Lemma}

\def \dim{{\mbox {dim}}\,}

\def\V{\mbox{Var}}

\def\Z{{\mathbb Z}}
\def\R\re
\def\V{\bf V}

\def \re{{\mathbb R}}

\def \0{\lambda_{0}}

\def\h{{\rm h}_{\rm top}(g)}

\begin{document}
\title[Growth rate of closed geodesics]{On the growth rate of contractible closed geodesics on
reducible manifolds}

\author[G. P. Paternain]{Gabriel P. Paternain}\thanks{G. P. Paternain
 was partially
 supported by CIMAT, Guanajuato, M\'exico}
 \address{ Department of Pure Mathematics and Mathematical Statistics,
University of Cambridge,
Cambridge CB3 0WB, England}
 \email {g.p.paternain@dpmms.cam.ac.uk}

\author[J. Petean]{Jimmy Petean}
 \address{CIMAT  \\
          A.P. 402, 36000 \\
          Guanajuato. Gto. \\
          M\'exico.}
\email{jimmy@cimat.mx}

\thanks{J. Petean is supported by grant 37558-E of CONACYT}


\date{May 2004}


\begin{abstract} We prove exponential growth rate of
contractible closed geodesics for an arbitrary bumpy metric on 
manifolds of the form $X_1\# X_{2}$ where 
the fundamental group of $X_1$ has a subgroup of finite index
at least 3 and $X_2$ is simply connected
and not a homotopy sphere.
\end{abstract}

\maketitle

\section{Introduction}

Let $M$ be a closed manifold and $g$ a smooth Riemannian metric on $M$.
Given $t\geq 0$, define $N_{g}(t)$ to be the number of geometrically
distinct closed geodesics of $g$ with length $\leq t$.
Similarly, let $N^0_{g}(t)$ be the number of geometrically
distinct {\it contractible} closed geodesics of $g$ with length $\leq t$.

If $g$ is a {\it bumpy metric}, i.e., all the closed geodesics of $g$ are 
non-degenerate, then $N_g(t)$ is finite
for each $t$ and we can ask the following basic question: what is the
behaviour of $N_g(t)$ as $t$ tends to infinity?
Similarly if $g$ is a {\it 0-bumpy metric}, i.e., all the {\it contractible}
closed geodesics of $g$ are non-degenerate, then $N^{0}_g(t)$ is finite
for each $t$ and we can also wonder about its growth.

Let $\Lambda(M)$ be the space of piecewise differentiable closed curves
$c:\re/\Z\to M$, endowed with the compact-open topology.

\medskip

\noindent {\bf Proposition A.} {\it Let $M$ be a closed manifold and
let $X\subset M$ be a simply connected submanifold, possibly with
boundary.
Consider  the inclusion map $\iota: X\to M$  and let $R_i$ be the
rank of the map induced in the loop space homology 
for some field of coefficients $k_p$, $p$ prime or zero.
Then, for any 0-bumpy metric $g$ on $M$ there exist constants 
$\alpha=\alpha(g)>0$ and $\beta=\beta(g)>0$ such that
\[N_g^0(t)\geq \alpha\, \max_{i\leq \beta t} R_i \]
for all $t$ sufficiently large.}

\medskip

This proposition is a generalization of a theorem of
M. Gromov \cite{G1}, later improved by W. Ballmann and W. Ziller \cite{BZ}, 
who proved the proposition when $M$ is simply connected and $X=M$
(we remark that in Gromov's theorem it is essential to assume that
$M$ is simply connected). In fact, the methods in \cite{PP2}
show that the proposition is still true if we replace $X$ by an arbitrary 
finite simply connected CW-complex $K$ and $\iota$ by any continuous map 
$f:K\to M$.

The result has the following interesting consequence:

\medskip

\noindent {\bf Theorem B.} {\it Let $M$ be a closed manifold of
  dimension $n\geq 3$. Suppose that $M$ can be 
decomposed as $X_1\#  X_2$, 
where $\pi_{1}(X_1)$ has a subgroup of finite index $\geq 3$ and $X_2$ 
is simply connected.
Then, for any 0-bumpy metric $g$ on $M$, $N_g^0(t)$ grows 
exponentially with $t$ unless $X_2$ is a homotopy 
sphere.}

\medskip

The novelty in Theorem B lies in the fact that $M$ is allowed to have
any fundamental group as long as it has subgroups of finite index $\geq 3$.
Results of this kind, when $M$ is simply connected were obtained
by P. Lambrechts in \cite{pascal} and in fact, we will rely quite heavily
on his Hochschild homology computations.

We note that there are finitely presented infinite groups without
proper subgroups of finite index. An example is given by the
Higman 4-group, see \cite{S}. We do not know if Theorem B is still
true for such groups.

A celebrated theorem of D. Gromoll and W. Meyer \cite{GM} asserts that
any Riemannian metric on a closed simply connected manifold $M$ has
infinitely many geometrically distinct closed geodesics if the Betti
numbers of $\Lambda(M)$ are unbounded for some field of coefficients.
Their methods combined with our proof of Theorem B yield the existence
of infinitely many geometrically distinct {\it contractible} closed geodesics
for {\it any} metric on a manifold as in Theorem B.
We do not know if the theorem is still true if we drop the bumpy condition
and we just assume $N^{0}_{g}(t)$ finite for all $t$.

There are several papers establishing lower bounds for the growth of
$N_g(t)$ in the presence of fundamental group \cite{Ball,BTZ,BaHi,BaK,H}.
These bounds are not exponential (e.g. $t/\log t$)
and do not give information about {\it contractible} closed geodesics, but
they do hold for {\it any} Riemannian metric.

The present note is a spin-off of our investigations on the topological entropy
$\h$ of the geodesic flow \cite{PP,PP1}. Except in the case of surfaces, there is basically no relationship between positivity of $\h$ and exponential growth 
for $N_g(t)$. M. Herman gave in \cite{He} an example of a minimal real
 analytic diffeomorphism on a closed 4-manifold with positive topological entropy.
However, if the geodesic flow of $g$ has a horseshoe, then
the number of {\it hyperbolic} closed geodesics of $g$ grows exponentially.
It has been conjectured that generically in the $C^k$ topology, $2\leq k\leq\infty$, the geodesic flow of $g$ has a horseshoe, but this has only been proved
for surfaces \cite{CP,KW}.

\medskip

{\it Acknowledgements}: We thank Wolfgang Ziller and
Werner Ballmann for various comments on the first draft of the manuscript.

\section{A brief review of Morse theory of the loop space}

Let $\Lambda=\Lambda(M)$ be the space of piecewise differentiable closed curves
$c:S^1=\re/\Z\to M$, endowed with the compact-open topology, i.e. the topology
induced by the metric $\rho(c,c')=\max_{t\in S^1}d(c(t),c'(t))$.
The energy functional $E:\Lambda\to\re$ is defined by $E(c)=\frac{1}{2}\int_{0}^{1}\langle \dot{c}(t),\dot{c}(t)\rangle\,dt$.
The critical points of $E$ are the closed geodesics and the point curves.
If $\ell(c)$ denotes the length of $c$, then the Cauchy-Schwarz inequality implies that $\ell(c)^2\leq 2\,E(c)$ with equality if and only if $c$ is parametrized proportional to arc length.

There is a natural action of $S^1$ on $\Lambda$. Given $s\in S^1$ and $c\in\Lambda$, we can define $sc$ by $sc(t):=c(t+s)$.
If $c$ does not reduce to a point, then the isotropy group of $c$ is a finite
subgroup of $S^1$ isomorphic to $\Z_k$ for some $k\geq 1$.
In this case, $k$ is called the {\it multiplicity} of $c$ and we can write
$c=d^k$, where $d^k(t)=d(kt)$. The curve $c$ is called the
 {\it $k$-th iterate} of $d$.

If $c$ is a closed geodesic, then the entire orbit of $c$ under the 
$S^1$-action is a set of critical points of $E$. A closed geodesic $c$ is said to be {\it non-degenerate} if the orbit of $c$ is a non-degenerate critical submanifold of $\Lambda$. Equivalently, $c$ is non-degenerate if and only if there is no periodic Jacobi field along $c$ orthogonal to $\dot{c}$.
In more dynamical terms, $c$ is non-degenerate if and only if the linearized Poincar\'e map of the orbit of the geodesic flow corresponding to $c$ does not have 1 as an eigenvalue.

A metric is said to be {\it bumpy} (resp. {\it 0-bumpy}) if all closed (resp. contractible)
geodesics are non-degenerate.
The bumpy metric theorem asserts that the set of $C^{r}$ bumpy metrics
is a residual subset of the set of all $C^{r}$ metrics endowed with
the $C^{r}$ topology for all $2\leq r\leq\infty$. The bumpy metric
theorem is traditionally attributed to R. Abraham \cite{Ab}, but see also
Anosov~\cite{An1} and Klingenberg and Takens \cite{KT}. 
Obviously bumpy implies 0-bumpy.

\section{Proof of Proposition A}

We will need the following version of a result due to 
M. Gromov \cite{G1} for manifolds with non-empty boundary,
see also \cite[p. 102]{pansu}, \cite{bates,Ba, P1}.
In all these references the manifold is assumed to have
empty boundary. In the case of the pointed loop space $\Omega (X)$
we have written down a proof of the theorem for compact
manifolds with non-empty boundary in \cite{PP1}.
For completeness we will show in the appendix that the 
details work as well for the free loop space $\Lambda (X)$.

\begin{Theorem}Given a metric $g$ on a simply connected
compact manifold $X$ (possibly with boundary), there exists a
constant $C=C(g)>0$ such that any element
in $H_{i}(\Lambda(X),k_p)$ can be represented by a cycle
whose image lies in $\Lambda^{Ci}(X)$. \label{coolfact}
\end{Theorem}

Let us prove now Proposition A. Let $\Lambda_0\subset\Lambda$ denote
the subspace of contractible closed loops. If all the contractible
closed geodesics are non-degenerate, then given $a>0$, there exist
only finitely many such closed geodesics $c_1,\dots,c_r$ with energy precisely
$a$. There also exists an $\varepsilon>0$ such that there is no contractible closed
geodesic $c$ with $E(c)\neq a$ but $a-\varepsilon\leq E(c)\leq a+\varepsilon$.
It follows from Morse theory that $\Lambda_0^{a+\varepsilon}$
is homotopy equivalent to $\Lambda_0^{a-\varepsilon}$ with cells $e_i$, $\bar{e}_{i}$
attached, where $\dim e_i=\dim \bar{e}_{i}-1=$ index of $c_i$.
It follows from this that $(\Lambda_0^a,M)$ is homotopy equivalent to a relative
CW-complex, where the number of $k$-cells attached to $M$ is equal to the number
of critical circles $S^1\,c$ of $E$ with index of $c$ equal to $k$ or $k-1$
and $0<E(c)\leq a$. Observe that the iterates of a closed geodesic
may give rise to different critical circles of $E$ of index $k$ or $k-1$.

Let $b_i(t):=\dim H_i(\Lambda_0^{t},M)$. It follows from the above that
the number of prime and iterate contractible closed geodesics of index $i$ or $i-1$
and length $\leq t$ is $\geq b_i(t)$.

We will need the following lemma proved in \cite{BZ} that gives control
of the contribution of a closed geodesic $c$ and its iterates to $b_i(t)$.
The lemma is not properly stated as such in \cite{BZ}, but the statement
can be found at the end of the proof of their main theorem.

\begin{Lemma} There exists a constant $n_0=n_0(g)$ such that at most
$n_0$ iterates of a contractible closed geodesic of length $\leq t$
can have index $i$ or $i-1$ if $t/2< i\leq t$.
\label{ballzi}
\end{Lemma}

From the lemma it follows that $N_g^0(t)\geq \max_{t/2<i\leq t}b_{i}(t)/n_{0}$.
Since $N_g^0(t)\geq N_g^0(t/2)$ we obtain
$N_g^0(t)\geq \max_{1<i\leq t}b_{i}(t)/n_{0}$. 

Proposition A is now a consequence of the following lemma:

\begin{Lemma} There exists a constant $C=C(g)>0$ such that
$b_i(t)\geq R_i$ for $n+1<i\leq C\,t.$
\end{Lemma}

\begin{proof}First note that $H_i(\Lambda_0^{t})\cong H_i(\Lambda^{t}_0,M)$ for $i>n+1$, so we
only need to prove that 
$b_i(\Lambda_0^{t})\geq R_i$ for $i\leq C\,t$.

Now take an element $e\in H_i(\Lambda(X))$. 
By Theorem \ref{coolfact} (we consider $X\subset M$ with the induced
metric) 
there exists $C>0$ (depending on $g$) such that whenever
$C\,t\geq i$, there exists a cycle $\eta\in \Lambda^{t}(X)$ such that 
$\eta$ represents
$e$. Since $\iota$ preserves energy, $\iota\circ\eta$ is a cycle
in $\Lambda_{0}^{t}$ and the lemma follows.

\end{proof}

\section{Proof of Theorem B}

\begin{proof} Let $\widetilde{X_1}$ be a finite covering of $X_1$
of degree $k\geq 3$. 
Consider the  covering $N$ of $M=X_1 \# X_2$ given by 
the connected sum of $\widetilde{X_1}$ with $k$ copies
of $X_2$. 
Since the exponential growth rate of closed geodesics is invariant
under finite coverings, we only need to work with $N$.
We now want to find a submanifold of $N$ to play the role of $X$
in Proposition A.
If $D$ is an $n$-dimensional disc then
$X:=D\# 3X_2$ appears as a submanifold of $N$.

\begin{Lemma} The inclusion map $\iota: X\to N$ induces a map
$$\iota_{*}:H_{*}(\Lambda(X),k_p)\to H_{*}(\Lambda(N),k_p)$$ 
whose rank grows exponentially for some $p$.
\label{topology}
\end{Lemma} 

\begin{proof} There is a map $p: M \# 3X_2 \to 3X_2$ obtained by
collapsing the exterior and the boundary of a disc $D\subset M$
(which contains the 3 points where the connected sums are taken)
to a point. If $j= p \circ \iota$ it is enough to prove that the
rank of the map induced by $j$ in the free loop space homology grows 
exponentially.

Clearly, $j$ is just the inclusion of $D\#3X_2$ in $3X_2 =S^n \#
3X_2$. In particular, the computations done in
\cite{pascal} apply to this problem. 
If $X_2$ is not a homotopy sphere, there is
a minimal degree in which its homology is non-trivial for some 
coefficient field $k_p$. In $3X_2$ we have then three
cohomology classes in this minimal degree which are independent and
whose cup product vanish. These cohomology classes imply the existence
of elements in the corresponding DGA's which satisfy the hypothesis
of \cite[Proposition 9]{pascal}. Of course, $j$ sends these cohomology
classes to corresponding ones in $3X_2 -D$. Then one looks at the map
induced in Hochschild homology and the  exponential
growth of the rank follows.

\end{proof}

\noindent
Theorem B follows directly from the previous lemma and Proposition
A.

\end{proof}

\section{Appendix: Proof of Theorem \ref{coolfact} when $\partial X$
 is non-empty}

\begin{proof}: We will only modify the proof in \cite{PP1} to show
that it works in the case of the free loop space.
Assume that the boundary of $X$ is non-empty. We have a collar
of the boundary diffeomorphic to $\partial X \times
[0,2)$. Let $Y$ be the manifold obtained by
deleting $\partial X \times [0,1)$ from $X$.
Of course, $Y$ is diffeomorphic to $X$. Now find a finite
number of convex subsets of $X$ which cover $Y$. Call
them $V_{\alpha}$, $1\leq \alpha  \leq k_0$. Let ${\bf T}$ be
a triangulation of $Y$. For a point $p\in Y$, let $F(p)$ be
the closed cell of lowest dimension containing $p$ and let
$O(p)$ be the union of all closed cells intersecting $F(p)$.
Note that $O(p)$ is a compact subset of $Y$.
Similarly, for any subset $K\subset Y$ one can define $F(K)$
as the union of $F(p)$ for $p\in K$ and $O(K)$ as the union
of $O(p)$ for $p\in K$. 
There is a positive number $\delta$ (depending on the covering
$\{ V_{\alpha} \}$ and $g$) such that, 
after taking some barycentric
subdivisions, we can take ${\bf T}$ so that for every subset $K$
of diameter bounded by $\delta$,   $O(K)$ is
contained in one of the $V_{\alpha}$'s.

Now we define the open subset $\Lambda_k$ of
$\Lambda (\cup V_{\alpha})$ as the set of all 
paths $\omega$ in $\cup V_{\alpha}$ with 
$\omega (0)=\omega (1)$ 
so that  for all $j$ between 1 and $2^k$,

$$O(\{ \omega(j-1/2^k ),\omega (j/2^k ) \} ) \ 
\cup \  \omega [j-1/2^k ,j/2^k ]$$

\noindent
is contained in one of the $V_{\alpha}$'s.
It is easy to see that $\Lambda (Y)$ is
contained in the union of the $\Lambda_k$'s. 

Let $B_k$ be the set of sequences $p_0 ,...,p_{2^k}$ of points 
in $Y$ such that $p_0 =p_{2^k}$ and for each $j$ between
1 and $2^k$ $O(\{ p_{j-1} ,p_j \} )$ is contained
in one of the $V_{\alpha}$'s. 
Let $\Lambda_k^Y \subset \Lambda_k$ be the set of paths in
$\Lambda_k$ for which all $2^{k}$ points $\omega(j/2^{k})$ are
in $Y$. Then $B_k$ is naturally identified with a
subset of $\Lambda_k^Y$ (an element of $B_k$ uniquely 
determines a broken geodesic which sends $j/2^k$ to $p_j$) and 
it is actually a deformation retract of $\Lambda_k^Y$.

Given a cycle representing a homology class in
$\Lambda (X)$ it can be represented by a cycle
in $\Lambda (Y)$ which is therefore contained in $\Lambda_k^Y$
for some $k$. Therefore we can retract it to $B_k$.
But $B_k$ is easily identified with a subset of
$Y^{2^k}$. Moreover, under this identification if a point
$(p_1 ,...,p_{2^k})\in B_k$, then the whole 
$F(p_1 )\times ...\times F(p_{2^k})$ is contained in
$B_k$. This implies that ${\bf T}$ induces a cell decomposition in
$B_k$. Hence the $i$-homology class can be represented by
a combination of cells of dimension $i$. A cell in $B_k$ is
a product of cells in each coordinate. The dimension of such
a cell is the sum of the dimensions of the corresponding
cells, of course. If the total dimension is $i$ then there can
be at most $i$ cells of positive dimension. 
Since $X$ is simply connected there exists a smooth map
$f:X\rightarrow X$ which is smoothly homotopic to the
identity and which sends the union of the images of all
the geodesic segments joining vertices 
in the triangulation to a point. The norm of the differential
of $f$ is bounded since $X$ is compact and 
$f$ induces a map $\hat{f} : \Lambda (X)\rightarrow \Lambda
(X)$ which is homotopic to the identity.
Now paths belonging to an $i$-cell of $B_k$
are formed by pieces joining vertices of the triangulation and
at most $2i$ pieces in which one of the points is not
a vertex. Under the map $\hat{f}$ the former are sent to a point
and the latter to a path of length bounded in terms of the
norm of the differential of $f$ and the diameter of the 
$V_{\alpha}$'s.
Therefore, there exists a constant $C(g)$
such that the image of the $i$-skeleton of
$B_k$ is sent by $\hat{f}$ to the subset of paths with
energy bounded by $C(g) \ i$. The theorem follows.

\end{proof}

\end{document}